\documentclass[11pt]{article}
\usepackage{amsmath,amsfonts,amssymb,amscd}
\usepackage{colordvi}   %%%  ADDED
\usepackage[bf,small,center,compact,topmarks,calcwidth,pagestyles]{titlesec} %%add new style

\def\N{{\mathbb N}}
\def\Z{{\mathbb Z}}

\def\l{\lambda}

\def\hpi{\widehat{\pi}}
\def\hK{\widehat K}

\def\mfS{{\mathfrak S}}
\def\x{{\mathbf x}}
\def\y{{\mathbf y}}

\def\moins{\raise 1pt\hbox{{$\scriptstyle -$}}}
\def\plus{\raise 1pt\hbox{{$\scriptstyle +$}} }
\def\s{\scriptstyle}
\def\bmoins{\fbox{\raise 1pt\hbox{{$\scriptstyle -$}}1}}
\def\demi{\frac{\s 1}{\s 2}}

\newcommand{\qed}{{\hfill\rule{4pt}{7pt}}}

\newtheorem{theorem}{Theorem}

\newtheorem{proposition}[theorem]{Proposition}
\newtheorem{lemma}[theorem]{Lemma}
\newtheorem{corollary}[theorem]{Corollary}
\newtheorem{definition}[theorem]{Definition}

\def\Proof{\noindent{\it Proof.}\ }

\begin{document}

\begin{center}
{\bf \large Non symmetric Cauchy kernels for the classical Groups}
\end{center}

\begin{center}
 Amy M. Fu$^*$\footnotetext{$^*$ supported by the 973 Project on Mathematical
Mechanization, the National Science Foundation, the Ministry of
Education, and the Ministry of Science and Technology of China. }
and Alain Lascoux
\end{center}

\begin{abstract}
We give non-symmetric versions of the Cauchy kernel and Littlewood's
kernels, corresponding to the types $A_n$, $B_n$, $C_n$ and $D_n$,
of the classical groups. We show that these new kernels are diagonal
in the basis of two families of key polynomials (one of them being
Demazure characters) obtained as images of dominant monomials under
isobaric divided differences. We define scalar products such that
the two families of key polynomials are adjoint to each other.
\end{abstract}

%\noindent{\bf Keywords:} Cauchy kernel; Littlewood's formulas;
%Weyl character formula; isobaric divided differences; key
%polynomials; Demazure characters.

%%%%%%%%%%%%%%%%%%%%%%%%%%%%%%%%%%%%%%
\section{Introduction}

Given two sets of indeterminates $\x =\{x_1,x_2,\ldots,x_n\}$,
$\y=\{y_1,y_2,\ldots,y_n\}$, the classical \emph{Cauchy kernel}
$\widetilde\Omega^{A}$ diagonalizes in the basis of Schur
functions~:
\begin{equation}          \label{Cauchy}
\widetilde\Omega^{A}=\prod_i\prod_j (1-x_i y_j)^{-1}=  \sum_{\l}
s_\l( \x)\, s_\l(\y).
\end{equation}

The Cauchy kernel may be considered as the generating function of
all characters of the symmetric groups. Multiplying the kernel
$\widetilde\Omega^{A}$ by the factor
\begin{equation*}
    \prod_{i<j} (1-x_ix_j) \prod_{i,j}(1-x_i/y_j)^{-1}  \quad \text{or} \quad
 \prod_{i \leq j} (1-x_ix_j)\prod_{i,j}(1-x_i/y_j)^{-1} ,
\end{equation*}
Littlewood \cite{Littlewood} obtained  expansions for the
following kernels
$\widetilde\Omega^{C}$ and $\widetilde\Omega^{D}$, in terms of
symplectic Schur functions and orthogonal Schur functions (see
below for the precise definitions):
 \begin{eqnarray}
   \widetilde\Omega^{C}=\frac{\prod_{i<j} (1-x_ix_j)}{\prod_i\prod_j
 (1-x_iy_j)(1-x_i/y_j)  }
 &=&
  \sum_{\l} s_\l(\x)\, Sp_\l(\y') ,   \label{Cauchy1}\\
  \widetilde\Omega^{D}=\frac{\prod_{i\leq j} (1-x_ix_j)}{\prod_i\prod_j
   (1-x_iy_j)(1-x_i/y_j)}
 &=&
  \sum_{\l} s_\l(\x)\, {\cal O}_\l(\y')  \label{Cauchy2} \, ,
\end{eqnarray}
where $\y'=\{y_1,y_2,\ldots,y_n,
y_1^{-1},y_2^{-1},\ldots,y_{n}^{-1}\}$.

In this paper, we shall study the following non-symmetric versions
of the kernels $\widetilde\Omega^{A}$, $\widetilde\Omega^{C}$ and
$\widetilde\Omega^{D}$~:
\begin{eqnarray*}
 \Omega^A &:=& \frac{1}{\prod_{i+j\leq n+1} ( 1-x_iy_j) },  \\
 [2pt]
 \Omega^B &:=& \frac{\prod_{i<j} (1-x_ix_j) \prod_i (1 + x_i)}
     {\prod_i\prod_j  (1-x_iy_j)  \ \prod_{i\leq j}(1-x_i/y_j)},  \\
     [2pt]
 \Omega^C &:=& \frac{\prod_{i<j} (1-x_ix_j) }{
     {\prod_i\prod_j  (1-x_iy_j)  \ \prod_{i\leq j}(1-x_i/y_j)}} ,  \\
     [2pt]
  \Omega^{D} &:=&\frac{\prod_{i\leq j} (1-x_ix_j) }{
     {\prod_i\prod_j  (1-x_iy_j)  \ \prod_{i\leq j}(1-x_i/y_j)}}.
\end{eqnarray*}

It will be convenient to interpolate between $\Omega^B$ and
$\Omega^C$, choosing an arbitrary parameter $\beta$, and defining~:
\begin{equation*}
 \Omega^{BC} =
  \frac{\prod_{i<j} (1-x_ix_j) \prod_i (1 +\beta x_i)}
     {\prod_i\prod_j  (1-x_iy_j)  \ \prod_{i\leq j}(1-x_i/y_j)}.
  \end{equation*}

For each type $A,B,C,D, BC$, there exist two families of isobaric
divided differences, which allow, starting from all dominant
monomials, to define two families of key polynomials, one of them
being the Demazure characters. Our main result (Th.
\ref{th:CauchyABCD}) is that all kernels $\Omega^A,\ldots,
\Omega^{BC}$ diagonalize in the corresponding basis of key
polynomials.

Notice that in type $A$, one also has a polynomial kernel, which is
the resultant $\prod_i\prod_j (x_i-y_j)$ of  two $z$-polynomials
$\prod_i (z-x_i)$ and $\prod_j (z-y_j)$. It still decomposes without
multiplicity in the basis of products of Schur functions in $\x$ and
$\y$. The non-symmetric version of the resultant, $\prod_{i,j:\,
i+j\leq n+1} (x_i-y_j)$, decomposes in the basis of products of
\emph{Schubert polynomials} in $\x$ and $\y$, and the main
properties of Schubert polynomials are easy consequences of the fact
that $\prod_{i,j:\, i+j\leq n+1} (x_i-y_j)$ is a reproducing kernel
\cite{Cbms}.

In the present article, for type $A$, we have rather taken  the
inverse function $\prod_{i,j:\, i+j\leq n+1} (1-x_iy_j)^{-1}$.  The
corresponding polynomials are no more the Schubert polynomials,
though there are interesting relationships between them and Demazure
characters.

The Cauchy kernel may be used to define a scalar product on the ring
of symmetric polynomials with coefficients in $\Z$, with respect to
which Schur functions constitute the only orthonormal basis
\cite{Macdonald}.  Starting from Weyl's denominators, we also define
scalar products with respect to which, for all types, the bases of
key polynomials are adjoint of each other (Th.
\ref{th:AdjointBases}). However, Bogdan Ion \cite{Ion01,Ion04} has
shown that key polynomials can be obtained as a limit case of
Macdonald polynomials. Thus the definition of the scalar product and
the orthogonality property of key polynomials result from the theory
of Macdonald polynomials. Nevertheless, we are giving an independent
derivation in sections 6, 7, because this approach relies only on
simple properties of divided differences and does not require double
affine Hecke algebras.

\section{Weyl Groups}

We shall realize the classical groups as groups operating on
vectors, or, equivalently, on Laurent polynomials, when considering
the vectors to be exponents of monomials. For more informations
about Coxeter groups, see \cite{BjornerBrenti}.

Fixing a positive integer $n$, we  define the operators $s_i$
$(1\leq i\leq n)$,  and $\tau_n$ acting on vectors $v \in \Z^n$ as
follows (\emph{operators are noted on the right}):
 \begin{eqnarray*}
v s_i &=& [\dots, v_{i+1},v_i,\ldots] \ , \ 1\leq i < n, \\
v s_n &=& [\dots, v_{n-1}, -v_n],   \\
v \tau_n &=& [\dots, -v_n,\, -v_{n-1}]   \, .
\end{eqnarray*}

Denoting a Laurent monomial $x_1^{v_1}\cdots x_n^{v_n}$ by $x^v$,
we extend by linearity the preceding operators to operators on
Laurent polynomials in indeterminates $x_1,\ldots, x_n$. The
simple transpositions $s_i$ $(i=1,\ldots, n-1)$ interchanges
$x_i$ and $x_{i+1}$, $s_n$ transforms $x_n$ into $x_n^{-1}$, and
$\tau_n$ sends $x_{n-1}$ onto $x_n^{-1}$,   $x_n$ onto
$x_{n-1}^{-1}$.

The group generated by $s_1,\ldots, s_{n-1}$ is a faithful
representation of
 the symmetric group $\mfS_n$ (type $A_{n-1}$). Adding the generator
$s_n$ gives the Weyl group of type $B_n$ or $C_n$ (which will be
distinguished later), while $s_1,\ldots, s_{n-1},\tau_n$ induce a faithful
representation of the type $D_n$.

An element $w$ of any of these groups can be identified with the
image under $w$ of the vector $v= [1,2,\ldots,n]$. For type
$A_{n-1}$, one gets \emph{permutations}; for type $B_n$, $C_n$, one
gets the \emph{bar-permutations}, writing $\bar r $ rather than
$-r$; and for  type $D_n$, one gets the \emph{bar-permutations with
an even number of bars}. As usual, we denote by $\ell(w)$ the
\emph{length} of $w$ (i.e. the length of a reduced decomposition of
$w$). For a given type, we extend the definition of length to
vectors: the length $\ell(v)$ of $v\in \Z^n$ is the minimum number
of generators of the group that must be applied to pass from $v$ to
the decreasing reordering of $|v_1|,\ldots, |v_n|$.

 There is a
unique element of maximal length for each type, usually denoted
$w_0$. For $A_{n-1}$, it
is $\omega^A:=[n,\ldots,1]$. For $B_n$, $C_n$, it is
$\omega^B=\omega^C:=[-1,\ldots, -n]$. For $D_n$, it is
$\omega^D:=[-1,\ldots, -n]$ if $n$ is even, and otherwise, it is
$\omega^D:=[-1,\ldots,-n+1, n]$. Reduced decompositions for these
elements are
\begin{eqnarray*}
  \omega^A &=& (s_1)\, (s_2 s_1)\, \cdots (s_{n-1}\cdots s_1)  \, , \\
 \omega^B &=&\omega^C= (s_n)\, (s_{n-1} s_n s_{n-1})\, \cdots
       (s_1\cdots s_{n-1} s_n s_{n-1} \cdots s_1)  \, ,   \\
 \omega^D &=& (s_{n-1} \tau_n)\, (s_{n-2} s_{n-1}\tau_n s_{n-2})\, \cdots
       (s_1\cdots s_{n-2}  s_{n-1}\tau_n s_{n-2} \cdots s_1)  \, .
\end{eqnarray*}

We shall also need  some conjugates $\theta_i$ of $s_n$
($i=1,2,\ldots,n$), defined by
$$  v \theta_i  = [\ldots, -v_i,\ldots]  \, . $$

In the group algebra, the most important element is the alternating
sum of all elements $ \sum_w   (-1)^{\ell(w)} w $. From the explicit
representation of the groups given above, it is easy to obtain the
following factorizations~:
\begin{equation}\label{BN}
 \sum_{w\in B_n, C_n} (-1)^{\ell(w)} w  =
 (1-\theta_1)\cdots (1-\theta_n) \,
 \sum_{\sigma\in \mfS_n} (-1)^{\ell(\sigma)} \sigma
\end{equation}
and
\begin{equation}\label{CN}
\sum_{w\in D_n} (-1)^{\ell(w)} w  = \frac{1}{2}
    \Big( (1-\theta_1)\cdots (1-\theta_n) + (1+\theta_1)\cdots
    (1+\theta_n)\Big)\, \sum_{\sigma\in\mfS_n} (-1)^{\ell(\sigma)} \sigma \, .
\end{equation}

\section{ The Weyl character formula}

In this section, we give a brief review of the Weyl character
formula, from an algebraic point of view only.

Let $\rho^A = \rho^D :=[n-1,\ldots, 1,0] $,
$\rho^B:=[n-\demi,\ldots,2-\demi,1-\demi]$, $\rho^C:=
[n,\ldots,2,1]$.  The sums
$$\sum_w  (-1)^{\ell(w)}  \left( x^{\rho^\heartsuit}   \right)^w\, , $$
$\heartsuit= A,B,C,D$, under the appropriate group, can be written
as determinants~:
\begin{eqnarray*}
 \Delta^A  &=& \det\Bigl( x_i^{j-1}\Bigr)_{1\leq i,j\leq n}\, ,  \\
  \Delta^B &=&
   \det \Big(x_{i}^{j-1/2} - x_{i}^{1/2-j}\Big)_{1\leq i,j\leq n} \, , \\
  \Delta^C &=& =\det \Big(x_{i}^{j} - x_{i}^{-j}\Big)_{1\leq i,j\leq n} \, , \\
  2 \Delta^D  &=& \det \Big(x_{i}^{j-1} + x_{i}^{-j+1} \Big)_{1\leq i,j \leq n}
                                                  \, .
\end{eqnarray*}

These determinants are easily factorized~:
\begin{eqnarray}
\Delta^A  &=& \prod_{i<j} (x_i-x_j) \, , \\
\label{DeltaB} \Delta^B &=&\prod_i \left(x_i^{1/2} -x_i^{-1/2}
\right)\,
                     \prod_{i<j} (x_i-x_j) (1-\frac{1}{x_i x_j})\, ,  \\
\label{DeltaC} \Delta^C &=&\prod_i \left(x_i -x_i^{-1}  \right)\,
                      \prod_{i<j} (x_i-x_j) (1-\frac{1}{x_i x_j}) \,,  \\
\label{DeltaD} \Delta^D &=& \prod_{i<j} (x_i-x_j) (1-\frac{1}{x_i
x_j})\, .
\end{eqnarray}

Taking now the images of general dominant monomials,
one obtains Weyl's expressions of the characters of the linear,
symplectic or orthogonal groups \cite{Weyl}~.
For $\l$ dominant, with $\ell(\lambda) \leq n$
in types $A,B,C$ and $\ell(\lambda)<n$ in type $D$, the quotient
$$ \left( \sum_w (-1)^{\ell(w)} x^{(\l+\rho)w} \right)
 \left( \sum_w (-1)^{\ell(w)} x^{\rho w} \right)^{-1}  $$
is equal to
\begin{eqnarray}                                \label{WeylA}
  s_\l(\x) &,&   \text{type } A  \, ,\\
                                                 \label{WeylC}
  Sp_\lambda(\x') &,&   \text{type } C \, ,\\
                                                 \label{WeylB}
  {\cal O}_\lambda(\x'') &,&   \text{type } B \, ,\\
                                                    \label{WeylD}
  {\cal O}_\lambda(\x') &,&   \text{type } D \, ,
\end{eqnarray}
where $\x'=\{x_1,\ldots,x_n, x_1^{-1},\ldots,x_{n}^{-1}\}$,
$\x''=\{x_1,\ldots,x_n, 1,x_1^{-1},\ldots,x_{n}^{-1}\}$.

For a combinatorial interpretation in terms of lattice paths, we
refer to Chen, Li, Louck \cite{CLL}.

In the remainder of this text, we shall be concerned with the generalization
of these characters by Demazure.

%%%%%%%%%%%%%%%%%%%%%%%%%%%%%%%%%%%%%%%%%%%%%%%%%%%%%%%%%%%%%%%%%%%%%%
\section{Divided differences and key polynomials}

Restricting to $n=1,2$, one can interpret Weyl's formulas as
operators on the ring of polynomials in one or two variables. These
operators are similar to Newton's divided differences. They are
called \emph{Demazure operators} \cite{Demazure}, or
\emph{isobaric divided differences}.

More specifically, for each type $A,B,C,D$, one defines two
families of divided differences acting on functions of
$x_1,\ldots, x_n$, and written on the right.

The first family is
\begin{eqnarray*}
\pi_i &:& \ f \longmapsto  f\, \pi_i :=\frac{x_if
-x_{i+1}f^{s_i}}{x_i-x_{i+1}} \ ,\, 1\leq i <n\, , \\
 \pi_n^C &:& \ f \longmapsto  f\, \pi_n^C :=
 \frac{x_n f -x_n^{-1} f^{s_n}}{x_n- x_n^{-1}} \, , \\
 [5pt]
 \pi_n^B &:& \ f \longmapsto  f\, \pi_n^B :=
 \frac{x_n f -f^{s_n}}{x_n- 1} \, , \\
 [5pt]
 \pi_n^D &:& \ f \longmapsto  f\, \pi_n^D :=
 \frac{ f - x_{n-1}^{-1}x_n^{-1}f^{\tau_n}}{1- x_{n-1}^{-1} x_n^{-1}} \, .
\end{eqnarray*}

 It is convenient to interpolate between  the operators $\pi_n^B$ and $\pi_n^C$
and define~:
$$ \pi_n^{BC} : \ f(x_1,x_2,\ldots,x_n) \longmapsto
\frac{(x_n+\beta) f -(x_n^{-1}+\beta) f^{s_n}}
  {x_n- x_n^{-1}} \, . $$
One sees that $\pi_n^C$ is recovered by putting $\beta=0$, while
$\pi_n^B$ corresponds to $\beta=1$. This operator results from the
 representation of the Hecke algebra of type
$\widetilde C_n$, defined by Noumi
(cf.  Sahi \cite[2.4]{Sahi}).

The second family is
$$   \hpi_i:= \pi_i -1 \, ,\ 1\leq i<n\, , $$
and
$$\hpi_n^\heartsuit= \pi_n^\heartsuit -1\ , \heartsuit=B,C,D,BC \, .$$

Each family satisfy the braid relations for type $A,B,C,D$
respectively \cite{Demazure}.
Notice that the operators $ \pi_i$ (resp. $\hpi_i$) , $1\leq i
\leq n$, commute with the multiplication by functions invariant
under $s_i$, and that $\pi_n^{D}$ (resp. $\hpi_n^D$)
 commutes with  the multiplication by functions invariant under
$\tau_n$. Thus, computations with a single $ \pi_i$, $i\leq n$ are
reduced to an action on the linear span of $1, x_i$.  In particular,
it is immediate to obtain that each operator satisfies the following
quadratic relations (which are degenerate cases of the Hecke
relations).

\begin{lemma} The squares of the isobaric divided differences satisfy
\begin{eqnarray*}
  \pi_i\pi_i = \pi_i\, ,  \  &&\ \hpi_i\hpi_i = -\hpi_i \, ,\, 1\leq i<n \, ,\\
  \pi_n^\heartsuit \pi_n^\heartsuit= \pi_n^\heartsuit\, ,   \ &&
   \ \hpi_n^\heartsuit \hpi_n^\heartsuit = -\hpi_n^\heartsuit   \, ,
   \heartsuit=B,C,D,BC\ .
\end{eqnarray*}
\end{lemma}

We define the \emph{key polynomials of type $\heartsuit$},  for
$\heartsuit= A,B,C,D,BC$, to be the images of dominant monomials
under products of isobaric divided differences. For type
$A,B,C,D$, these are the \emph{Demazure characters}. Using the
divided differences $\hpi_i$ instead of $\pi_i$, one obtains a
second family of key polynomials.

In more details. The starting points of all families are
$$ x^\l =x_1^{\l_1}x_2^{\l_2}\ldots x_n^{\l_n}
 = K_\l^\heartsuit = \hK_\l^\heartsuit\, ,
\quad \text{all partitions } \l\in \N^n \, .$$ The other polynomials
are defined recursively by
\begin{equation}
  K_v^\heartsuit\, \pi_i
            = K_{v\, s_i}^\heartsuit\, \ \& \
  \hK_v^\heartsuit\, \hpi_i
      = \hK_{v\, s_i}^\heartsuit\, ,\,\text{when}\ v_i>v_{i+1} \, ,\, i<n\, .
\end{equation}
\begin{equation}
  K_v^\heartsuit\, \pi_n^\heartsuit
                     =     K_{v\, s_n}^\heartsuit\,  \ \&\
  \hK_v^\heartsuit\, \hpi_n^\heartsuit
           = \hK_{v\, s_n}^\heartsuit\, ,\,\text{when}\ v_n>0 \, ,\,
  \text{for}\ \heartsuit=B,C,BC\,  .
\end{equation}
\begin{equation}
  K_v^D\, \pi_n^D
         = K_{v\, \tau_n}^D  \ \&\
  \hK_v^D\, \hpi_n^D
        = \hK_{v\, \tau_n}^D\, ,\, l(v\tau_n)>l(v).
\end{equation}

The definition is consistent since the operators satisfy the braid
relations. Notice that, when $v\in \N^n$, then all
$K_v^\heartsuit$ (resp. $\hK_v^\heartsuit$), $\heartsuit=
A,B,C,D,BC$ coincide with each other,
 since the exceptional generators are not used in the computation.
 In that case, we
shall write $K_v$, $\hK_v$, ignoring the types. We shall also
need to use at the same time operators acting on $x_1,\ldots,
x_n$, and operators acting on $y_1,\ldots, y_n$. In that case, we
use superscripts.

The images of a dominant monomial $x^\lambda$ under the maximal
divided difference $\pi_\omega^\heartsuit$, for $\heartsuit=
A,C,B,D$, are respectively the RHS of Eq. (\ref{WeylA}),
(\ref{WeylC}), (\ref{WeylB}), (\ref{WeylD}).

For $\heartsuit=BC$, and $\beta=-1$, one recovers
the \emph{odd symplectic characters} of Proctor \cite[Prop. 7.3]{Proctor}.

Divided differences can be extended to operators on paths. We refer
specially to the work of Littelmann \cite{Littelmann1,Littelmann2}.

%%%%%%%%%%%%%%%%%%%%%%%%%%%%%%%%%%%%%%%%%%%%%%%%%%%%%%%%%%%%%%%%%%%%
\section{Cauchy-type Kernels}

In this section, we shall show that all the kernels $
\Omega^{\heartsuit}$, $\heartsuit= A,B,C,D,BC$, are diagonal in the
basis of key polynomials. In fact, our computations will essentially
be reduced to the following cases, the verifications of which are
immediate.
\begin{eqnarray}
\label{Tool5} (1- ax_i)^{-1}\, \pi_i
                &=&  (1- ax_i)^{-1} (1- ax_{i+1})^{-1} \\[3pt]
\label{Tool6} (1- ax_i)^{-1}\, \hpi_i
                &=&  a x_{i+1} (1- ax_i)^{-1} (1- ax_{i+1})^{-1} \\[3pt]
 \label{Tool1}
 (1- ax_{i+1})\, \pi_i
      &=&  (1-a /x_i)\, \pi_i \ = \ 1 \, ,\ 1\leq i<n\, ,  \\[3pt]
  \label{Tool2}
 (1- ax_{i+1}) (1-b /x_i)\, \pi_i &=&   1- ab  \, , \ 1\leq i<n\, , \\[3pt]
  \label{Tool3}
 (1- b/x_n)\, \pi_n^{BC}  &=& 1+ \beta b \, , \\[3pt]
  \label{Tool4}
  (1-b /x_{n-1})(1-b /x_n)\, \pi^D_n      &=&
  1-b^2 .
\end{eqnarray}

%To compute the image of rational functions, we shall get rid of
%denominators as follows. Given $i$, and $ f/g$ a rational function,
%$f$ being a rational function invariant under $s_i$, and $g$ a
%polynomial, then
%\begin{equation}{\label{RA}}
%\frac{f}{g}\, \pi_i = g^{s_i}\pi_i\,\frac{f}{ g g^{s_i}} \, .
%\end{equation}

We introduce the operator
$$ \Xi_n := \sum_{\sigma\in\mfS_n}
                \hpi_\sigma^x \, \pi_{\sigma\omega}^y   \, ,$$
where $\omega$ is the maximal element in $\mfS_n$. Filtering the
set of permutations according to the position of $n$,
 one gets the following factorization.
\begin{lemma}
 We have
\begin{equation}
  \Xi_n = \Xi_{n-1}\, \left( \sum_{i=0}^{n-1}
   \hpi_{[n\moins1 : i]}^x\,  \pi^y_{[n\moins1 : n\moins 1\moins i]} \right)\,
   ,
\end{equation}
where
 $$ \pi_{[n\moins 1:i]} := \pi_{n-1}\, \pi_{n-2} \cdots \pi_{n-i} \, .$$
\end{lemma}

For example, the element  $\Xi_4$ factorizes as
 $$ \Xi_4 = \Xi_3\, \left(
   \pi_3^y\pi_2^y\pi_1^y +   \hpi_3^x \pi_3^y\pi_2^y  +
    \hpi_3^x \hpi_2^x \pi_3^y  + \hpi_3^x \hpi_2^x \hpi_1^x \right) \, .$$

The next proposition shows that the operator $\Xi_n$ allows to obtain
the kernel $\Omega^A$ from the generating function of the dominant
monomials.

\begin{proposition} {\label{pro:Xi}} We have
\begin{multline}
\frac{1}{(1-x_1y_1)   (1-x_1x_2y_1y_2) \cdots
  (1-x_1 \cdots x_n y_1\cdots y_n)}\ \Xi_n    \\
  = \frac{1}{\prod_{i+j\leq n+1}  1-x_iy_j }\ =\ \Omega^A  \, .\
\end{multline}
\end{proposition}

\Proof  The factor $(1-x_1 \cdots x_n y_1\cdots y_n)^{-1}$
commutes with all the divided differences $\hpi_i^x$, $\pi_i^y$,
$1\leq i\leq n-1$.
 Using the above factorization of $\Xi_n$, and supposing the proposition
true for $n-1$, one has to compute the image of $\prod_{i+j\leq
n}(1-x_iy_j)^{-1}$  under the sum
$$\sum_{i=0}^{n-1} \hpi_{[n\moins1 : i]}^x\,
\pi^y_{[n\moins1 : n\moins 1 \moins  i]}.$$

By repeated use of (\ref{Tool6}), one obtains
$$ \prod_{i+j\leq n}
(1-x_iy_j)^{-1} \, \hpi_{n-1}^x \cdots  \hpi_k^x =
  \prod_{i+j\leq n}
(1-x_iy_j)^{-1} \frac{ x_ny_1}{1-x_ny_1} \cdots
  \frac{ x_{k+1} y_{n-k}  }{1- x_{k+1} y_{n-k}}\, .  $$
Thanks to (\ref{Tool5}), the action of $\pi_{n-1}^y \cdots
\pi_{n-k+1}^y $ on this last function reduces to multiplication by
$$  \frac{1}{1-x_1y_n}\, \frac{1}{1-x_2 y_{n-1}} \cdots
    \frac{1}{1-x_{k-1} y_{n-k+2}}   \, .$$

Reducing now the  sum to a common denominator,
 it can be rewritten as  the
product of  $\Omega^A$  times the factor
$$ \sum_{k=1}^{n-1}x_n\ldots x_{k+1}y_1\ldots
y_{n-k}(1-x_ky_{n-k+1})+(1-x_ny_1).
$$
This last factor is nothing but the factor $ (1-x_1 \cdots
x_n y_1\cdots y_n)$ which commutes with all the divided differences.
This completes the proof.  \qed

\begin{lemma}{\label{Le1}} Let $\Phi^{BC}$  be the following
 operator acting on the variables $y_1, \ldots,y_n$.
$$ \Phi^{BC} :=  \bigl(\pi_n^{BC} \pi_{n-1} \cdots \pi_1\bigr)\,
   \bigl(\pi_n^{BC} \pi_{n-1} \cdots \pi_2\bigr)
 \cdots
  \bigl(\pi_n^{BC} \pi_{n-1}\bigr)\, \bigl(\pi_n^{BC} \bigr) .
$$
Then
$$ \Omega^A \, \Phi^{BC}  = \Omega^{BC}.  $$
\end{lemma}

\Proof Each step of the computation of $ \Omega^A \, \Phi^{BC}$
corresponding to the above factorization of $\Phi^{BC}$, is of one
of the following two types.

When $i<n$, then the rational function can be written
$(1-x'y_i)^{-1}(1-x/y_{i+1})^{-1}\, f$, with $f=f^{s_i}$, i.e. $f$
symmetrical in $y_i$ and $y_{i+1}$.  In that case, thanks to
(\ref{Tool2}), one has
\begin{eqnarray*} \frac{f}{(1-x'y_i)(1-x/y_{i+1})}
\pi_i &=&
 \frac{(1-x'y_{i+1})(1-x/y_i) \pi_i \, f }{
 (1-x'y_i) (1-x'y_{i+1})(1-x/y_{i+1}) (1-x/y_i)} \\
 &=& \frac{(1-xx')\, f}{
 (1-x'y_i) (1-x'y_{i+1})(1-x/y_{i+1})(1-x/y_i) } \, .
 \end{eqnarray*}
In other words, this step has consisted in multiplying  by
$$  (1-xx') \, (1-x'y_{i+1})^{-1}(1-x/y_i)^{-1}  \, .$$

In the case of a step $\pi_n^{BC}$, the current function is
$(1-xy_n)^{-1}\, f$, where $f=f^{s_n}$. Thanks to (\ref{Tool3}),
one has
\begin{multline*}
 \frac{f}{1-xy_n}\, \pi^{BC}_n  =
   (1-x/y_n)\, \pi^{BC}_n\, \frac{f}{(1-xy_n) (1-x/y_n)}  \\
 =  (1 + \beta x) \, \frac{f}{(1-xy_n) (1-x/y_n)}  \,
\end{multline*}
and therefore, we have created the factor $  (1 + \beta x)\,
(1-x/y_n)^{-1}$.

The product of all the factors we have created is
 $$ \prod_i (1 + \beta x_i)\,
 \prod_{i<j} (1-x_ix_j)\prod_{i+j>n+1} (1-x_iy_j)^{-1}\prod_{i\leq j}(1-x_i/y_j)^{-1} \, , $$
which is indeed equal to the quotient to $\Omega^{BC}/ \Omega^A$.
\qed

To treat the type $D$, we  define recursively  the following
operators (still acting on $y_1,\ldots, y_n$ only)~:
\begin{multline}
 \Phi^D_2 = \pi_1 \pi^D_2 ,\ \Phi^D_3 = \left( \pi_2 \pi^D_3\right)\,
  \pi_1 \pi_2 \pi^D_3\, , \ldots, \\
  \Phi^D_n  = \left( \Phi^D_{n-1} \right)^+ \,
  \pi_1 \pi_2\cdots \pi_{n-1} \pi^D_n  \, ,
\end{multline}
where the symbol $(\, )^+$ denotes the shift $i\to i+1$ of all
indices inside the parentheses. For example, taking $n=4$, we have
$$  \Phi^D_4 = \left(\pi_3 \pi_4^D\right)\,
  \left(\pi_2\pi_3 \pi_4^D\right)\, \left(\pi_1\pi_2\pi_3 \pi_4^D\right)\, .$$

\begin{lemma} {\label{Le2}}
Let
$$
\Omega^A_{n-1}= \prod_{i+j\leq n} (1-x_iy_j)^{-1}.
$$
Then
$$  \Omega^A_{n-1} \, \Phi^D_n  = \Omega^D    \, .$$
\end{lemma}

\Proof   The successive steps in the computation of the image of
$\Omega^A_{n-1}$ are of three possible types.

Step $\pi_{n-1}$. The current function is $(1-xy_{n-1})^{-1} f$,
with $f$ is symmetrical in $y_n$ and $y_{n-1}$.  Thanks to
(\ref{Tool5}), we have
$$ \frac{f}{1-xy_{n-1}} \pi_{n-1} =
\frac{f }{(1-xy_n)(1-xy_{n-1})}\, .
$$
We have just created a factor $ (1-xy_n)^{-1} $.

Step $\pi_n^{D}$. The current function is $(1-xy_{n-1})^{-1}
(1-xy_n)^{-1} f$, with $f$ invariant under $\tau_n$. Thanks to
(\ref{Tool4}),   we have
\begin{multline*}
 \frac{f}{(1-xy_{n-1})(1-xy_n)}\, \pi_n^{D}
  =  \frac{(1-x /y_{n-1})(1-x /y_n)\pi_n^{D} \, f}
  {(1-xy_{n-1})(1-xy_n)(1-x /y_{n-1})(1-x / y_n)}\\
 =  \frac{ (1-x^2)\, f}{(1-xy_{n-1})(1-xy_n)(1-x /y_{n-1})(1-x / y_n)}  \,
\end{multline*}
and the transformation is just multiplication by
$$  (1-x^2)\, (1-x /y_{n-1})^{-1}(1-x /y_n)^{-1} \, .$$

Step $\pi_i$, $i<n-1$. The current function is
$(1-x'y_i)^{-1}(1-x/y_{i+1})^{-1}\, f$, with $f=f^{s_i}$. Thanks
to (\ref{Tool2}), one has
\begin{multline*}
\frac{f}{(1-x'y_i)(1-x/y_{i+1})}\, \pi_i = \frac{
(1-x'y_{i+1})(1-x/y_i) \pi_i\, f}{
 (1-x'y_i) (1-x'y_{i+1})(1-x/y_i) (1-x/y_{i+1})} \\
 = \frac{ (1-xx')\, f}{
 (1-x'y_i) (1-x'y_{i+1})(1-x/y_i) (1-x/y_{i+1})} \, .
\end{multline*}
The function has been multiplied by
$$  (1-xx') \, (1-x'y_{i+1})^{-1}(1-x/y_i)^{-1}  \, .$$

The products of all  the above factors is
 $$ \prod_{i\leq j} (1-x_ix_j)
     \prod_i\prod_j (1-x_iy_j)^{-1} \prod_{i\leq j}(1-x_i/y_j)^{-1}\, . $$
This completes the proof. \qed

The preceding relations between the different kernels, and the
function $(1-x_1y_1)^{-1}   (1-x_1x_2y_1y_2)^{-1}\cdots$ allow to
expand these kernels.
\begin{theorem}   \label{th:CauchyABCD}
We have
\begin{eqnarray}  \label{CauchyA}
  \Omega^A &=&   \sum_{v\in \N^n} \hK_v(\x)\, K_{v\omega}(\y),    \\
 \label{CauchyBC}
  \Omega^{BC} &=& \sum_{v\in \N^n} \hK_v(\x)\, K_{-v}^{BC}(\y),    \\
 \label{CauchyD}
  \Omega^D &=& \sum_{v\in \N^n:\, v_n=0} \hK_v(\x)\, K_{-v}^{D}(\y),
\end{eqnarray}
where $x_n$ is specialized to $0$ in the last equation.
\end{theorem}

\Proof Note that
$$\frac{1}{(1-x_1y_1)   (1-x_1x_2y_1y_2) \cdots
  (1-x_1 \cdots x_n y_1\cdots y_n)}$$ is the generating function
of all dominant monomials $x^\l y^\l$ in $n$ indeterminates
$x_1y_1,\ldots, x_ny_n$.

From the definition of $\Xi_n$ and of key polynomials, one has
$$  \sum_\l  x^\l y^\l\, \Xi_n =
  \sum_{v\in \N^n} \hK_v(\x)\, K_{v\omega}(\y) \,,  $$
where the sum ranges over all partitions $\l$ of length at most
$n$. Thus,  Proposition \ref{pro:Xi}, entails (\ref{CauchyA}).

The image of a key polynomial $K_{v_n, \ldots, v_1} (\y)$ , $v\in
\N^n$  under $\Phi^{BC}$ is $K_{-v}(\y)$. Therefore, the image of
the RHS of (\ref{CauchyA}) under $\Phi^{BC}$ is the RHS of
(\ref{CauchyBC}), and Lemma \ref{Le1} gives (\ref{CauchyBC}).

Similarly, the image of $K^y_{v_{n-1},\ldots, v_1, 0}$ under
$\Phi^D_n $ is $K^y_{-v_1,-v_2,\ldots, -v_{n-1}, 0}$. Therefore,
the image of the expansion of $\Omega^A_{n-1}$ under $\Phi^D_n $
is the RHS of (\ref{CauchyD}), and Lemma (\ref{Le2}) completes the
proof of (\ref{CauchyD}) and of the theorem.    \qed

Note that (\ref{CauchyA}) has been established  combinatorially
  in \cite{DoubleCrystal},
using  Schensted bijection and double crystal graphs.

Let us conclude this section by showing  that the identities
(\ref{CauchyA}), (\ref{CauchyBC}) and (\ref{CauchyD}) imply the
Cauchy formula and Littlewood's formulas respectively.
Indeed, $\hpi_i \pi_i=0$, $1\leq i <n$,
a fortiori,  $\hpi_i \pi_\omega=0$, where $\omega$ is the maximal
element of $\mfS_n$.
Therefore all the summands in the right hand sides of
(\ref{CauchyA}), (\ref{CauchyBC}) and (\ref{CauchyD}) are sent to
$0$ under $\pi_\omega^x$,  except the terms
$$  \hK_\l(\x)\, \pi_\omega^x = K_{\l\omega}(\x) =
 s_\l(\x).  $$

On the other hand,
\begin{eqnarray*}
\prod_{i+j\leq n+1} (1-x_iy_j)^{-1}\, \pi_{\omega}^x  &=&
 \prod_{i,j=1}^n (1-x_iy_j)^{-1}  \\
\prod_{1\leq i\leq j \leq n} (1-x_i/y_j)^{-1}\, \pi_{\omega}^x  &=&
 \prod_{i,j=1}^n  (1-x_i/y_j)^{-1}
\end{eqnarray*}

 Specializing $\beta$, we get Cauchy formula (\ref{Cauchy}) and Littlewood's identities
 (\ref{Cauchy1}), (\ref{Cauchy2}), as images of
(\ref{CauchyA}), (\ref{CauchyBC}) and (\ref{CauchyD})
respectively.

\section{Scalar products}

Bogdan Ion \cite{Ion01, Ion04} has shown how to obtain the two
families of Demazure characters $K_v^\heartsuit$,
$\hK_v^\heartsuit$, $\heartsuit=A,B,C,D$, by degeneration of
Macdonald polynomials. Degenerating also Cherednik's scalar product
\cite{Cherednik}, one gets a  scalar product for each of the types
$\heartsuit$,  with respect to which the bases $\{K_v^\heartsuit\}$,
$\{\hK_v^\heartsuit\}$ are adjoint of each other. But instead of
having recourse to the elaborate theory of non symmetric Macdonald
polynomials, we shall directly define scalar products on
polynomials, and check orthogonality properties by simple
recursions.

Recall that in the theory of Schubert polynomials \cite{Cbms},
one defines a scalar product by using the maximal divided difference;
as a consequence divided differences are self-adjoint.
 This scalar product can also be written
$$ (f,g) = CT\left( fg
            \prod_{1\leq i ,j\leq n} (x_i^{-1} -x_j^{-1}) \right)\, ,$$
where $CT$ means ``constant term''.

It is easy to adapt this definition to our present needs,
keeping the compatibility of the scalar product with the isobaric
divided differences.

One first replaces the Vandermonde determinant by Weyl's
denominators, $\Delta^B$, $\Delta^C$ and $\Delta^D$ multiplied by
$x^{\rho}$. We add to their list
\begin{equation}  \label{DeltaBC}
 \Delta^{BC}:= \Delta^C \prod\nolimits_{i=1}^n (1+\beta x_i)^{-1} \, ,
\end{equation}
keeping  $\rho^{BC}=\rho^C =[n,\ldots,1]$.

\begin{definition}
 For $\heartsuit=B,C,D,BC$, and for  Laurent polynomials $f,g$ in
 $x_1,\ldots, x_n$, let
\begin{eqnarray}      \label{ScalProd}
 (f,g)^\heartsuit  &=&
    CT\left( fg \, x^{\rho^{\heartsuit}}\, \Delta^\heartsuit \right)  \, ,\\
    [5pt]             \label{ScalProdA}
(f,g)^A  &=&
   CT\Big( f(x_1,\ldots, x_n) g(x_n^{-1},\ldots,x_1^{-1})
                \prod_{1\leq i<j\leq n} (1-x_ix_j^{-1}) \Big),
\end{eqnarray}
where one expands $(1+\beta x_i)^{-1}$
 as a formal series in the variable $ x_i \beta$.
\end{definition}

For example, taking $n=2$, one has
\begin{eqnarray*}
(f,g)^A &=& CT\left( f(x_1,x_2)\, g(x_2^{-1},x_1^{-1}) \, (1-x_1
x_2^{-1})
  \right),\\
[8pt] (f,g)^{BC} &=& CT\left( \frac{fg\, x_1^2x_2\,
   (x_1-x_1^{-1}) (x_2-x_2^{-1}) (x_1-x_2)(1-x_1^{-1}x_2^{-1})}
 {  (1+ x_1\beta)(1+ x_2\beta)}\right),\\
[5pt] (f,g)^D &=& CT\Big( fg\, x_1\,(x_1-x_2)(1-\frac{1}{x_1x_2})
     \Big).
\end{eqnarray*}

Notice that the scalar product $(f,g)^{BC}$ does specialize to
$(f,g)^B$ for $\beta=1$, and to $(f,g)^C$ for $\beta=0$:
\begin{eqnarray*}
  (f,g)^B &=& CT \Big( fg x_1^{3/2}x_2^{1/2}(x_1^{1/2}-x_1^{-1/2})
  (x_2^{1/2}-x_2^{-1/2})(x_1-x_2)(1-\frac{1}{x_1x_2})\Big)\\
  &=& CT\Big( fg\, x_1 (x_1-1)(x_2-1)(x_1-x_2)(1-\frac{1}{x_1x_2})
  \Big),\\
  [5pt]
(f,g)^C &=& CT\left( fg\, x_1^2x_2\,
   (x_1-\frac{1}{x_1}) (x_2-\frac{1}{x_2}) (x_1-x_2)(1-\frac{1}{x_1x_2})
   \right).
\end{eqnarray*}

Weyl defined a scalar product on symmetric functions by:
$$   (f,g) = (n!)^{-1}\, CT\bigl(
 f(x_1,\ldots, x_n)\, g(x_1^{-1},\ldots, x_n^{-1})\Delta^2   \bigr)
\, ,$$
and similarly for the other types, taking the squares of $\Delta^\heartsuit$.
Thus, as for the different Cauchy kernels (\ref{Cauchy}),
(\ref{Cauchy1}), (\ref{Cauchy2}), one passes from the symmetric case
to the non-symmetric one by taking ``half'' of the factors.

The crucial property of the scalar products (\ref{ScalProd}) and
(\ref{ScalProdA}) is the following compatibility with isobaric
divided differences.

\begin{theorem} \label{CompatibilityPi} Write $\pi_n=\pi_n^{\heartsuit}$,
$\hpi_n=\hpi_n^{\heartsuit}$, for $\heartsuit=B,C,BC,D$. Then the
operators $\pi_i$ and $\hpi_i$ ($1\leq i \leq n$) are self-adjoint
with respect to $(\, ,\, )^\heartsuit$, i.e. for every pair of
Laurent polynomials $f,g$, one has
$$
  \bigl( f \pi_i \, ,\, g \bigr)^\heartsuit   =
   \bigl( f \, ,\, g\pi_i  \bigr)^\heartsuit \,,\quad
   \bigl( f \hpi_i \, ,\, g \bigr)^\heartsuit=
   \bigl( f \, ,\, g\hpi_i  \bigr)^\heartsuit \,.
$$

In the case of type $A$, for $1\leq i \leq n-1$, $\pi_i$ (resp.
$\hpi_i$) is adjoint to $\pi_{n-i}$ (resp. $\hpi_{n-i}$), i.e. for
every pair of Laurent polynomials $f,h$, one has
$$
\bigl( f \pi_i \, ,\, h \bigr)^A=\bigl( f \, ,\, h \pi_{n-i}
\bigr)^A, \quad \bigl( f \hpi_i \, ,\, h \bigr)^A=\bigl( f \, ,\,
h \hpi_{n-i} \bigr)^A.
$$
\end{theorem}

\Proof To treat  all types in a uniform way, we write
$h(x_n^{-1},\ldots,x_1^{-1}) = g(x_1,\ldots, x_n)$. Then $h
\pi_{n-i} = g\pi_i$, $1\leq i<n$, and
$$ (f, h \pi_{n-i})^A = CT \Big( f\, (g \pi_i) \prod_{1\leq
             i<j\leq n} (1-x_ix_j^{-1})\Big)  \, ,$$

For all types, and $i<n$, the scalar product can now be written as
 $$   CT\Bigl(
   CT_{x_i,x_{i+1}} \left(fg(1-x_i/x_{i+1})  \clubsuit     \right)  \Bigr),   $$
where $\clubsuit$ is a function symmetrical in $x_i,x_{i+1}$ and $
CT_{x_i,x_{i+1}}$  is the constant term in the variables
$x_i,x_{i+1}$ only.

Let us write $f,g$ as $f= f_1 +x_{i+1} f_2$, $g= g_1 +x_{i+1}
g_2$, with $f_1,f_2,g_1,g_2$ invariant under $s_i$. The difference
$ f\pi_i g - g\pi_i f= f\hpi_i g - g\hpi_i f$ is equal to $
(f_1g_2-g_1f_2) x_{i+1}$. Therefore the constant term
\begin{eqnarray*}
&&CT_{x_i,x_{i+1}} \Bigl(  (f\pi_i g - g\pi_i f )\,
                   (1-x_i/x_{i+1})\clubsuit \Bigr) \\
&& \qquad \qquad \qquad= CT_{x_i,x_{i+1}} \Bigl(  (f\hpi_i g -
g\hpi_i f )\,
(1-x_i/x_{i+1})\clubsuit \Bigr)\\
&&\qquad \qquad \qquad=  CT_{x_i,x_{i+1}} \Bigl( (x_i-x_{i+1})\,
(f_1g_2-g_1f_2) \clubsuit \Bigr)
\end{eqnarray*}
 is null, because the function inside parentheses
is antisymmetrical in $x_i,x_{i+1}$.

In the case $i=n$,  $\heartsuit=BC$, one writes
$$ (f \, ,\, g )^{BC} = CT\Big(
   CT_{x_n} \Big(fg\, \frac{x_n}{1 +\beta x_n} (x_n-{x_n}^{-1})
       \clubsuit     \Big)  \Big)\, ,   $$
where $\clubsuit$ is a function invariant under $s_n$. Therefore,
to evaluate $ (f \hpi_n \, ,\, g )^{BC} - (f \, ,\, g\hpi_n )^{BC}
    = (f \pi_n \, ,\, g )^{BC} - (f \, ,\, g\pi_n )^{BC} $,
one can first compute
$$ CT_{x_n} \Big(
 \left(f \hpi_ng -g\hpi_n f \right)
 \frac{x_n}{1 +\beta x_n} (x_n-x_n^{-1}) \clubsuit \Big)
 = CT_{x_n} \Big( (g^{s_n}f  - f^{s_n}g) \clubsuit \Big)  \,
$$
which is null, because the function under parentheses is
alternating under $s_n$.

Similarly,  for $\heartsuit=D$, neglecting a function invariant
under $\tau_n$, to determine $ (f \hpi_n \, ,\, g )^D - (f \, ,\,
g\hpi_n )^D = (f \pi_n \, ,\, g )^D - (f \, ,\, g\pi_n )^D $, one
can first compute
$$
  CT_{x_{n-1},x_n} \Big(( f\hpi_n g- g\hpi_n f) (1-x_{n-1}x_n)
   \Big)  = CT_{x_{n-1},x_n} \Big(  f^{\tau_n}g -g^{\tau_n}f
   \Big) \,
$$
which is also null, because the function $f^{\tau_n}g
-g^{\tau_n}f$ is
 alternating under $\tau_n$.
This completes the proof. \qed

\section{Orthogonality}

Let us  extend the usual dominance order on partitions
\cite{Macdonald}
 to an order on vectors in $\Z^n$. Given two vectors
$u=(u_1,u_2,\ldots,u_n)$ and $v=(v_1,v_2,\ldots,v_n)$ in $\Z^n$,
$u\leq v$
% or $u-v\leq 0$
means the following inequalities
$$ \quad u_1\leq v_1,\, u_1+u_2\leq v_1+v_2,\,
  u_1+u_2+u_3\leq v_1+v_2+v_3,\, \ldots .  $$
One also extends the notation $|\l|$ to vectors: $|v|:=v_1+\cdots
+v_n$.

We give in the following lemmas some easy properties of the scalar
product.

\begin{lemma}\label{Lemma2}
For every partition $\l$ in $\N^n$, for every element $w$ in the
groups of types $B_n$, $C_n$ and $D_n$ (resp. type
$A_{n-1}$), every monomial $x^u$  appearing in the expansion of
$x^\l\, \pi_w$ is such that $u\geq -\l$ (resp. $u \geq \l
\omega$).
\end{lemma}
\Proof  By recursion on length, one sees that $K^\heartsuit_v$,
$\heartsuit \neq A$, is equal to $x^v +\sum c_u^v x^u$, with
$v<u$.  \qed

\begin{lemma}\label{lemma1} For $u, v \in \Z^n$, $\heartsuit \neq A$,
\begin{equation*}
 (x^v , x^u)^\heartsuit \neq 0 \quad \text{implies that}
  \quad  v\leq -u  \, .
\end{equation*}
For $v, u \in \N^n$,
\begin{equation*}
(x^v , x^u)^A \neq 0  \quad \text{implies that} \quad
   v \leq u\omega, \quad \text{and} \quad |v|=|u| \, .
\end{equation*}
\end{lemma}

\Proof  Rewrite $x^\rho \Delta^C $ as the determinant
$$
\pm\det\Big(x^{j-i}(x_i^{2n-2j+2}-1)\Big)_{i,j=1}^n.
$$
If one expands the determinant by rows, then the powers of $x_1$
are nonnegative, the term $x_2^{-1}$ is multiplied by strictly
positive powers of $x_1$, the term $x_3^{-2}$ is mutliplied by
monomials in $x_1,x_2$ of degree at least $2$, $\cdots$.
Therefore, the scalar product $(x^v,x^u)^C$  can have  a constant
term only if
$$
v_1+u_1 \leq 0, v_1+v_2+u_1+u_2\leq 0, v_1+v_2+v_3+u_1+u_2+u_3\leq
0,\ldots,
$$
i.e. $v\leq -u$.

For $(x^v,x^u)^{\heartsuit}$, $\heartsuit=B,D$, the proof is
similar.
 For $(x^v,x^u)^{BC}$, we have multiplied $(x^v,x^u)^C$ by
formal series in $x_1,x_2,\ldots$ with positive exponents. Therefore
$(x^v,x^u)^{BC} \neq 0$ still implies $v\leq -u$.

For $(x^v,x^u)^A$, rewrite
 the product $\prod_{1\leq i<j \leq n}(1-x_ix_j^{-1})$ as the determinant
$\det(x_i^{j-i})_{i,j=1}^n$. One obtains that the scalar product
$(x^v,x^u)^A$  has a constant term only if
$$
v_1-u_n \leq 0, v_1+v_2-u_{n}-u_{n-1}\leq 0,\ldots,
$$
i.e. $v\leq u\omega$. Moreover, to have a non-zero constant term,
the total degree must be $0$, i.e. $|v|=|u|$. \qed

\begin{lemma}\label{Coro1}
For two partitions $\l$ and $\mu$ in $\N^n$, for $\heartsuit=A, B,
C, D, BC$, if there exists $w$ such that
$$ (x^\l\, \pi_w\, ,\, x^\mu)^\heartsuit \neq 0 \, ,$$
then $\l=\mu$.
\end{lemma}
\Proof There exists at least one monomial $x^v$ in $x^\l\, \pi_w$
such that $(x^v, x^\mu)^\heartsuit \neq 0$. Lemma \ref{lemma1}
then implies that $ v\leq -\mu$, and therefore $-\lambda \leq -\mu
$  when $\heartsuit\neq A$. In the case of type $A$, one has $v
\leq \mu \omega$ and therefore $\l \omega \leq \mu \omega$.

Write now $ (x^\l\, \pi_w\, ,\, x^\mu)^\heartsuit= (x^\l, x^\mu
\pi_{w'})^\heartsuit$, where  $w'$ is obtained from $w$ by using
Th. \ref{CompatibilityPi}. We have reversed the role of $\lambda$
and $\mu$, and therefore non-nullity of the scalar product implies
that $\l=\mu$. \qed

Note that if all the parts of $\lambda$ are different, then $w$
must be the maximal element of the group.

\begin{corollary} \label{coro2}  Let $\l$ be a partition in $\N^n$,
with at least one $0$ component when $\heartsuit=D$ and $n$ is
odd. Then
\begin{equation}\label{eq1}
 \heartsuit \neq A\, ,\   (K_v\, ,\, x^\l)^\heartsuit
         \neq 0 \quad \text{implies that}\ v=-\l  \, .
\end{equation}
In that case  $(K_{-\l}\, ,\, x^\l)^\heartsuit =1$.

\begin{equation}\label{eq12}
   (K_v\, ,\, x^\l)^A \neq 0 \quad \text{implies that}\ v=\l \omega  \, .
\end{equation}
In that case  $(K_{\l \omega}\, ,\, x^\l)^A =1$.
\end{corollary}

\begin{lemma}\label{LBCD}  Let $\heartsuit \neq A$ and $i:\, 1\leq i\leq n$.
Given four polynomials $f_1$, $f_2=f_1\pi_i$, $g_1$, $g_2= g_1
\hpi_i$, then $ (f_1,g_1)^{\heartsuit}=0 \ \& \
(f_2,g_1)^{\heartsuit}=1$, implies that
$$ (f_1,g_2)^{\heartsuit}=1 \qquad \& \qquad (f_2,g_2)^{\heartsuit}=0   \, . $$
Moreover,  any space $V$ stable under $\pi_i$ which is orthogonal
to $g_1$ is orthogonal to $g_2$.
\end{lemma}

\Proof We have
$$(f_2,g_2)^{\heartsuit}= (f_1 \pi_i\, ,\, g_1\hpi_i)^{\heartsuit}=
 (f_1 \pi_i\hpi_i, g_1)^{\heartsuit}= 0 \, ,$$
 and
 $$
  (f_1,g_2)^{\heartsuit}= (f_1,\, g_1\hpi_i)^{\heartsuit}
  = (f_1 , g_1\pi_i-g_1)^{\heartsuit} =(f_1\pi_i,g_1)^{\heartsuit}
  =(f_2, g_1)^{\heartsuit}=1 \, .$$

The last statement is immediate.  \qed

The next lemma has a similar proof.

\begin{lemma}\label{LA}
Given an integer $i$ $(1\leq i \leq n-1)$ and four polynomials
$f_1$, $f_2=f_1\pi_{n-i}$, $g_1$, $g_2= g_1 \hpi_i$, then $
(f_1,g_1)^A=0 \ \& \ (f_2,g_1)^A=1$ implies that
$$ (f_1,g_2)^A=1 \qquad \& \qquad    (f_2,g_2)^A=0   \, . $$
Moreover,  any space $V$ stable under $\pi_{n-i}$ which is
orthogonal to $g_1$ is orthogonal to $g_2$.
\end{lemma}

We are now ready to conclude.

\begin{theorem}  \label{th:AdjointBases}
 Let $u,v\in \Z^n$, and $\heartsuit
 \neq A$. Suppose moreover that, when $\heartsuit=D$ and
$n$ is odd, then $u$ or $v$ has at least one component equal to
$0$. Then
\begin{equation}\label{Neq1}
  (K_v\, ,\, \hK_u)^\heartsuit  = \delta_{-v,u}   \, ,
\end{equation}
where, as usual,  $\delta_{-v,u}$ is the Kronecker delta.

In the case of type $A$, for $u,v\in \N^n$, we have
\begin{equation}\label{Neq2}
(K_v\, ,\, \hK_u)^A  = \delta_{v\omega,u}   \, .
\end{equation}
\end{theorem}

\Proof  When $u$ is dominant, (\ref{eq1}) implies that $ (K_v\,
,\, x^{u})^\heartsuit=\delta_{-v,u}$. By induction on length,
suppose that $u$ is such that $\hK_u$ is orthogonal to every
$K_v$, except $(K_{-u}, \hK_u)^{\heartsuit}=1$. Let $i$ be such
that the linear span of $\hK_u, \hK_u\hpi_i$ be two-dimensional.
Then one uses Lemma \ref{LBCD},
 with $f_1= K_{-us_i}$,
$f_2=K_{-u}=f_1\pi_i$,  $g_1= \hK_u$, $g_2=\hK_{us_i}=g_1 \hpi_i$,
and $V$ generated by all $K_v$, $v\neq -u, -us_i$, to conclude
that $\hK_{us_i}$ is orthogonal to all $K_v$, except for $(
K_{-us_i}, \hK_{us_i})^{\heartsuit} =1$.

For type $A$, one replaces Lemma \ref{LBCD}
by  Lemma \ref{LA} to arrive to a similar
conclusion. \qed .

\vspace{0.3cm}

\noindent{\bf Acknowledgments.}
We are grateful to Bogdan Ion for his detailed explanations
about the connections between Macdonald polynomials and
Demazure characters.

\begin{flushleft}
\quad Center for Combinatorics, LPMC\\
\quad Nankai University, Tianjin 300071, P. R. China\\
\quad \emph{Email address:} fu@nankai.edu.cn

\vskip 0.5cm

\quad CNRS, IGM Universit\'e  de Marne-la-Vall\'ee\\
\quad 77454 Marne-la-Vall\'ee Cedex, France\\
\quad \emph{Email address:} Alain.Lascoux@univ-mlv.fr
\end{flushleft}£¬

\end{document}